\documentclass{amsart}

\usepackage{epsfig}

\newtheorem{thm}{Theorem}[section]
\newtheorem{lemma}[thm]{Lemma}

\newcounter{Tlistc}

\newcounter{ticklistc}

\newcommand{\ccc}{\theta}
\newcommand{\CCC}{\Theta}
\newcommand{\Int}{\text{Int}}

\newcommand{\Ima}{\text{Im}}
\newcommand{\Eul}{\text{Eul}}
\newcommand{\Tors}{\text{Tors}}

\newcommand{\cut}{\text{cut}}

\newcommand{\spn}{\text{spn}}

\newcommand{\aug}{\text{aug}}

\newcommand{\sing}{\text{sing}}
\newcommand{\Ker}{\text{Ker}}
\newcommand{\Rr}{\mathbb R}
\newcommand{\Qq}{\mathbb Q}
\newcommand{\Zz}{\mathbb Z}

\newcommand{\Hom}{\mathrm{Hom}}

\begin{document}

\title{A function on the homology of 3-manifolds}

\author{Vladimir Turaev}
\address{IRMA, CNRS et Universit\'e Louis Pasteur,
7 rue Ren\'e Descartes, 67084 Strasbourg Cedex, France}

    \address{Department of Mathematics, Indiana University,
    Bloomington, IN 47405, USA}

\date{\today}


\begin{abstract} In analogy with the Thurston norm, we define
for an  orientable 3-manifold $M$    a numerical function  on
$H_2(M;\Qq/\Zz)$. This function measures the minimal complexity of
folded surfaces representing a given homology class. A similar
function is defined on the torsion subgroup of $H_1(M;\Zz)$. These
functions are estimated  from below in terms of abelian torsions of
$M$.

\end{abstract}

\maketitle

\section{Introduction}

 One of the most   beautiful invariants of a
 3-dimensional manifold $M$ is the Thurston semi-norm on
 $H_2(M;\Qq)$, see \cite{th}.
The geometric idea leading to this semi-norm is to consider
 the minimal genus of a surface in $M$ realizing any given
2-homology class of $M$.
 Thurston's definition of the semi-norm uses a suitably
 normalized Euler characteristic of the surface rather than the genus. The
Thurston
 semi-norm   is uninteresting for a
 rational homology sphere $M$, since  then  $H_2(M;\Qq)=0$. However, a rational
 homology sphere  may have non-trivial 2-homology  with coefficients in   $\Qq/\Zz$.
Homology classes
  in $H_2(M;\Qq/\Zz)$    can be   realized by
 folded surfaces,   locally
 looking like unions of  several  half-planes in $\Rr^3$ with  common
 boundary line.  It is natural to consider
 \lq\lq smallest"  folded surfaces in  a given homology class.

 We use this train of ideas to define
for   an arbitrary
  orientable 3-manifold $M$ (not necessarily a rational homology
  sphere)
  a  function
 $$\ccc :H_2(M;\Qq/\Zz)\to \Rr_+=\{r\in \Rr\,\vert \, r\geq 0\}. $$
 This function measures the \lq\lq minimal" normalized Euler characteristic  of a  folded surface
 representing a given   class in $H_2(M;\Qq/\Zz)$.

 Using the
 boundary homomorphism $$d:H_2(M;\Qq/\Zz)\to H_1(M)=H_1(M;\Zz),$$  whose image is equal to $\Tors\, H_1(M)$, we derive from
 $\ccc$ a function $$\CCC:\Tors\, H_1(M)\to \Rr_+$$ by ${\CCC}(u)=\inf_{x\in d^{-1}(u)} \ccc(x)$
 for any $u\in \Tors\, H_1(M)$. One can view   $\CCC(u)$
    as a \lq \lq
  normalized
minimal genus" of oriented knots in $M$  representing   $u$. If $M$
is a rational homology sphere, then $d$ is an isomorphism and
$\CCC=\ccc\circ d^{-1}$.

 We give an estimate  of the function $\ccc $ from above in terms of the Thurston semi-norm
 on knot complements in $M$. This estimate implies that
 $\ccc $ is bounded from above and is upper semi-continuous  with respect to
 a natural
 topology on  $H_2(M;\Qq/\Zz)$.  (I do not know whether $\ccc $ is
 continuous.) The functions  $\ccc $ and $\CCC$ are also estimated
 from below  using abelian torsions of $M$.
 These estimates are parallel to the
  McMullen \cite{mc} estimate  of the Thurston semi-norm in terms of the
  Alexander polynomial.

 In contrast to the Thurston semi-norm, the   function $\ccc $ is
 non-homogeneous, that is in general $\ccc(kx)\neq k \,\ccc(x)$ for $k\in \Zz$
and
 $x\in H_2(M;\Qq/\Zz)$.
 Examples show that  the function $\ccc $ may not  satisfy  the triangle inequality.

The Thurston semi-norm of a 3-manifold $M$ is fully determined by
the Heegaard-Floer homology of $M$, see  \cite{oz}, and by the
Seiberg-Witten monopole homology of $M$, see \cite{km}. It would be
interesting to obtain similar computations of the functions $\ccc$
and $\CCC$.

The organization of the paper is as follows. We introduce the
functions $\ccc $ and $\CCC$ in  Section \ref{section:lk} and
estimate them  from above in Section \ref{section:2k}. In Section
\ref{section:3k} these functions are estimated from below in the
case
  where the first Betti number of the 3-manifold is non-zero. A similar
  estimate  for rational homology spheres is given in Section
  \ref{section:4k}. In Section 6 we describe a few examples. In Section
  7 we make several miscellaneous remarks.

Throughout the paper, the unspecified group of coefficients in
homology is   $\Zz$.

\section{Folded surfaces and the functions $\ccc ,\CCC$}\label{section:lk}

 \subsection {Folded surfaces}   By a
{\it  folded  surface} (without boundary), we   mean   a compact
2-dimensional polyhedron
  such that  each  point   has a neighborhood homeomorphic to a union  of
  several
       half-planes in $\Rr^3$ with  common
 boundary line. Such a neighborhood  is homeomorphic to     $
\Rr \times \Gamma_n $  where $n $ is a positive integer and
$\Gamma_n$ is a union of $n$ closed intervals with one common
endpoint  and no other common points.

The  {\it interior} $\Int (X)$ of a folded surface $X$ consists of
the points of $X$ which have   neighborhoods  homeomorphic  to
$\Rr^2$. Clearly, $\Int (X)$ is a 2-dimensional manifold.    The
{\it singular set} $\sing (X)=X- \Int (X)$ of $X$ consists of a
finite number of disjoint circles.  A neighborhood of a component
  of $\sing (X)$ in $X$ fibers over this component  with fiber $\Gamma_n$ for
some $n\neq 2$.

Cutting  out $X$ along $\sing (X)$ we obtain a   compact 2-manifold
 (with boundary) $X_{\cut}$. Each component   of $\Int (X)$ is the
interior of a  component    of $X_{\cut}$. Set
  $\chi_- (X)   =\sum_Y \chi_- (  Y)  $, where $Y$ runs over all components
of $X_{\cut}$
 and $$\chi_- (  Y)  = \max (-\chi  ( Y),0).$$   The
number $\chi_- (X)\geq 0$ measures the complexity of $X$. It is
 equal to zero if and only if all components of
 $X_{\cut}$ are either spheres or tori or annuli or disks.

 By an {\it orientation} of a  folded surface  $X$, we mean an orientation of the 2-manifold  $\Int (X)$.
 An  orientation of $X$ allows us to view $X$
  as a singular 2-chain with integer coefficients. This 2-chain is
   denoted   by the same letter  $  X $.
 Its boundary     expands as $ \sum_K \,i(K)\, \langle K\rangle$
 where $K$ runs over   connected components of $\sing (X)$, the
 symbol $\langle K\rangle$ denotes a 1-cycle on $K$ representing a
 generator of $H_1(K)\cong \Zz$, and  $i(K)\in \Zz$. Multiplying, if
 necessary, both $\langle K\rangle$ and $i(K)$ by $-1$, we can assume
 that $i(K) \geq 0$.
 In this way the   integer $i(K) $ is uniquely determined by $K$.
 It is called the \emph{index} of $K$ in $X$.
 For     $K$ with $i (K)\neq 0$, the 1-cycle $\langle K\rangle$
   determines an orientation of $K$.
   We say that this orientation is {\it induced} by the one on $X$.

We call a folded surface $X$ {\it simple} if it is oriented, the set
$ \sing (X)$ is homeomorphic to a circle, and its index in $X$ is
non-zero. This index is denoted   $i_X$. Note that $X$ is not
required to be connected; however, all components of $X$ but one are
closed oriented 2-manifolds.

\subsection {Representation of 2-homology by folded surfaces}
Let $M$ be an orientable 3-manifold. By a   folded surface in  $M$,
we mean a folded surface {\it embedded} in $ M$. Given a  simple
folded surface $X$  in  $M$,  the   2-chain
  $ (i_X)^{-1} X$   with rational coefficients
  is a 2-cycle modulo $\Zz$.    This cycle
   represents   a homology class in $H_2(M;\Qq/\Zz)$
   denoted~$[ X ]$.

The short exact sequence of groups of coefficients $0\to \Zz \to \Qq
\to \Qq/\Zz\to 0$ induces an exact homology sequence
\begin{equation}\label{seq}\cdots \to H_2(M;\Qq) \to
H_2(M;\Qq/\Zz)\to H_1(M )\to H_1(M;\Qq)\to \cdots.\end{equation}
 The homomorphism $  H_2(M;\Qq/\Zz) \to
H_1(M)$   in this sequence will be denoted $d_M$ and called the {\it
boundary homomorphism}.  The   exactness of  \eqref{seq} implies
that the image of $d_M $ is equal to the  group  $\Tors \, H_1(M) $
consisting of all elements of $H_1(M)$ of finite order.

For a  simple folded surface $X$  in $M$, the homomorphism $d_M$
sends $[ X ] $  into the 1-homology class represented by the circle
$\sing(X)$ with orientation induced by the one on $X$.

For example, if $X\subset M$ is a compact oriented 2-manifold with
connected non-void boundary, then $X$ is a simple folded surface
with $\sing (X)=\partial X$, $i_X=1$, and  $[X]=0$. Another example:
consider an unknotted circle $K$ lying in a 3-ball in $M$ and pick
$n\neq 2$ closed 2-disks bounded by $K$ in this ball and having no
other common points. We orient these disks so that the induced
orientations on $K$ are the same. The union  of these disks, $X=X(n)
$, is  a simple folded surface with $\sing (X)=K$, $i_X=n$, and $[ X
]=0$.

\begin{lemma}\label{le1}    Any homology class  $x\in H_2(M;\Qq/\Zz)$ can be
represented by a simple folded surface.
\end{lemma}

\begin{proof} Set $d=d_M:  H_2(M;\Qq/\Zz) \to H_1(M)$. We can  represent   $d(x)\in \Tors\, H_1(M)$ by an oriented embedded
circle $K\subset \Int (M)= M-\partial M$.    Pick an integer $n\geq
1$ such that $ n \, d(x)=0$. The standard arguments, using the
Poincar\'e duality and transversality, show that there is a simple
folded surface $X$  in $M $ such that $\sing (X)=K$ and $i_X=n$.

Since both $X$ and $M$ are orientable, the 1-dimensional normal
bundle of $\Int(X)$ in $M$ is trivial. Keeping $\sing(X)$ and
pushing $X-\sing(X)$ in a normal direction, we obtain a \lq\lq
parallel" copy $X_1$ of $X$ such that $X\cap
X_1=\sing(X_1)=\sing(X)=K$. The orientation of $X$ induces an
orientation of $X_1$ in the obvious way. Repeating this process
$k\geq 1$   times, we can obtain $k$ parallel copies $X_1, X_2,...,
X_k $ of $X$   meeting each other exactly at  $K$. Then $X^{(k)}=X_1
\cup X_2 \cup ... \cup X_k$ is a  simple folded surface such that
$\sing(X^{(k)})=K$ and  $i_{X^{(k)}}=   nk$. It follows from the
construction  that $[ X^{(k)} ]=[ X ]\in H_2(M;\Qq/\Zz)$ for all
$k\geq 1$.

The equalities $ d(x)=[K]=d ([ X ])$ imply that   $x- [ X ] \in \Ker
\, d=\Ima\, j$, where $j$ is  the   homomorphism $H_2(M;\Qq )\to
H_2(M;\Qq/\Zz)$ induced by the projection $\Qq\to \Qq/\Zz$. Pick
$y\in j^{-1} (x- [ X ]) \subset H_2(M;\Qq )$. There is an
  integer $k\geq 1$ such that   $ ky$ lies in the image of
the coefficient homomorphism $H_2(M;\Zz )\to H_2(M;\Qq )$. {\it A
fortiori}, the homology class $nky$ lies in this image.  We
represent $nky$ by a closed oriented (possibly non-connected)
surface $\Sigma \subset M$. Since $d(x)\in \Tors\, H_1(M)$, the
intersection number   $\Sigma \cdot K =\Sigma \cdot d(x)  $ is $0$.
Applying if necessary
  surgeries of index 1 to $\Sigma$, we can assume that $\Sigma
\cap K=\emptyset$.   Then $y$ is represented by the 2-cycle $
(nk)^{-1}\Sigma $ in $M-K$ and  $x=[ X ]+j(y)=[ X^{(k)} ] +j(y)$ is
represented by the 2-cycle $(nk)^{-1} (    X^{(k)}  + \Sigma ) $ mod
$\Zz$. Applying to $X^{(k)}$ and $\Sigma$ the usual cut and paste
technique, we can transform their union into a simple folded surface
$Z$ such that
  $\sing(Z)=\sing (X^{(k)})=K$ and $i_{Z}=nk$. Clearly,
$[ Z ]=x$.
\end{proof}

\subsection {Functions $\ccc $ and ${\CCC}$}\label{der} For an orientable 3-dimensional
manifold $M$, we define a function $\ccc =\ccc_M:H_2(M;\Qq/\Zz)\to
\Rr_+$ by
\begin{equation}\label{eq1} \ccc(x)=\inf_X \,  \frac
{\chi_- (X)}{i_X},\end{equation} where $x\in H_2(M;\Qq/\Zz)$ and $X$
runs over all simple folded surfaces
 in $M$    representing~$x$.
  In particular, the class $x=0$ can be
represented by the simple folded surface $X=X(n)\subset M$ with
$n\neq 2$, constructed before
  Lemma \ref{le1}. The equality  $\chi_- (X)=0$ implies that  $\ccc(0)=0$.

 For a
simple folded surface $X$,
   denote by $-X$ the same simple folded surface with opposite
   orientation in its interior.
The obvious equalities $$[ -X ]= -[ X ] ,\quad \chi_- (-X)=\chi_-
(X),\quad i_{-X}=i_X$$  imply that $\ccc(-x)=\ccc(x)$ for all $x\in
H_2(M;\Qq/\Zz)$.

  We  define a
function ${\CCC}={\CCC}_M: \Tors \, H_1(M)  \to \Rr_+$ by
\begin{equation}\label{eq2}{\CCC}(u)=\inf_{x\in d^{-1}(u)} \ccc(x) =\inf_X \,  \frac
{\chi_- (X)}{i_X},\end{equation} where $u\in \Tors \, H_1(M)$, $X$
runs over all simple folded surfaces in $M$ such that the circle
$\sing(X)$ represents $u$, and $d:H_2(M;\Qq/\Zz) \to H_1(M)$ is the
boundary homomorphism. In \eqref{eq2}, we can restrict ourselves to
connected $X$. Indeed, all components of $X$ disjoint from $\sing
(X)$ are closed oriented
  surfaces. They may be removed from $X$ without increasing
$\chi_- (X)$.

The properties of $\ccc$ imply that $\CCC(0)=0$ and
$\CCC(-u)=\ccc(u)$ for all $u\in \Tors \, H_1(M)$.  By the very
definition of $\CCC$, for all $x\in H_2(M;\Qq/\Zz)$,
$$\ccc(x)\geq \CCC(d(x)).$$

Using folded surfaces with boundary, we can similarly define
relative versions
$$H_2(M,\partial M;\Qq/\Zz)\to \Rr_+\quad {\text {and}}\quad \Tors\,
H_1(M,\partial M; \Zz)\to \Rr_+$$ of the functions $\ccc$ and
$\CCC$. We will not study them in this paper.

\subsection {Constructions and examples}\label{coae}
1. Let $\Sigma$ be a   closed connected 2-manifold embedded in an
oriented 3-manifold $M$. Let $K\subset \Sigma$ be a simple closed
curve such that   $\Sigma -K$ has an orientation which switches to
the opposite when one crosses $K$ in $\Sigma$. (Such an orientation
exists when $\Sigma$ is orientable and  $K$ splits $\Sigma$ into two
surfaces or when $\Sigma$ is non-orientable and $K$   represents the
Stiefel-Whitney class $w^1(\Sigma)\in H^1(\Sigma;\Zz/2\Zz)=
H_1(\Sigma;\Zz/2\Zz)$.) The  orientations of $M$ and $\Sigma-K$
induce an orientation of the normal bundle of $\Sigma -K$ in $M$.
Keeping $K$ and pushing $\Sigma-K$ in the corresponding  normal
direction, we obtain a   copy $\Sigma'$ of $\Sigma$ such that
$\Sigma'$ transversely meets $\Sigma$ along $K$.
  The union $X=\Sigma \cup \Sigma'$ is a simple folded
surface such that $\sing (X)= K$ and $i_X=4$.    Then  $\ccc([ X
])\leq (1/4)\, \chi_- (X)=  (1/2)\, \chi_- (\Sigma -K)$.

For example, we can apply this construction to the  projective plane
$\Sigma=\Rr P^2$ in $\Rr P^3$ taking as $K$ a
  projective circle on $\Rr P^2$. The resulting simple folded
surface $X$ represents the only non-zero  element $x$ of $H_2(\Rr
P^3;\Qq/\Zz) =\Zz/2\Zz$ because $\sing (X) $ represents the non-zero
element  of $H_1(\Rr P^3)=\Zz/2\Zz$. The equality $\chi_-(\Sigma
-K)=0$ implies  that $\ccc_{\Rr P^3}=0 $  and ${\CCC}_{\Rr P^3}=0 $.

2. Consider the 3-dimensional lens space $M=L(p,q)$, where $p,q$ are
co-prime integers with $p\geq 2$. The manifold $M$ splits as a union
of  two solid  tori  with common  boundary. It is easy to exhibit a
folded surface $X\subset M $
  such that $\sing (X)$ is the core circle of one of the solid tori  and
  $X-\sing (X)$ is a disjoint union  of $p$ open 2-disks. This
implies that the function ${\CCC}_M $ annihilates the elements of
$H_1(M)  $ represented by the core circles of the solid tori. Under
an appropriate isomorphism $H_1(M)\cong \Zz/p\Zz$, these   elements
correspond to $1\, ({\text {mod}} \, p)$ and $ q\, ({\text {mod}} \,
p)$. This implies that $\CCC_M=0$ if $p=2$ or $p=3$ or $p=5, q=2$.
For $p=2$, we recover the previous example, since $L(2,1)=\Rr P^3$.

3. Let  $K $ be  an oriented  homologically trivial knot in an
oriented 3-manifold $N$. Let $M$ be obtained by a $(p,q)$-surgery on
$K $
 where
 $p,q$ are co-prime integers with $p\geq 2$. Thus, $M$ is obtained by
cutting out a tubular neighborhood $U\subset N$ of $K$   and gluing
it back along a homeomorphism $\partial U\to \partial U$ mapping the
meridian $\mu\subset \partial U$ of $K$ onto a curve on $\partial U$
homological to $p \mu+ q \lambda $, where $\lambda\subset \partial
U$  is the longitude of $K$ homologically trivial in $N-K$.    The
element $u\in H_1(M)$ represented by the (oriented) core circle   of
the solid torus $U\subset M$ has finite order. This follows from the
fact that the $p$-th power of the core circle  is homotopic in
$U\subset M$ to    $\lambda \subset
\partial U$.  We claim that ${\CCC}(u)=  0$ if $K$ is a trivial knot
in $N$ and ${\CCC}(u)\leq p^{-1} (2g-1)$ if $K$ is a non-trivial
knot of  genus $g\geq 1$. Indeed, the longitude $\lambda$ bounds in
  $N- \Int (U)$ an embedded compact connected oriented surface of genus
  $g$.
This surface  extends in the obvious way to a simple folded surface
$X$ in $ M$ such that $\sing(X)$ is the core circle of $U\subset M$
and $i_X=p$. Clearly, $\chi_-(X)= \max(2g-1,0) $. This implies our
claim. (For $p=2$, one should \lq\lq double" $X$ along $\sing (X)$
as in Example 1.) As we shall see below,  if $K$ is a non-trivial
fibred knot and $p\geq 4g-2$, then ${\CCC}(u)= p^{-1} (2g-1)$.

\section{Estimates from above and semi-continuity}\label{section:2k}

 In this section we  estimate the function $\ccc=\ccc_M $ from above using  the Thurston norm.
Throughout this section, $M$ is a   connected orientable 3-manifold
(possibly, non-compact).

 \subsection{Comparison with the Thurston norm}
 Recall
first the  definition of the Thurston semi-norm $\Vert \cdot
\Vert_M$ on $H_2(M;\Qq)$.   The Poincar\'e duality (applied to
compact submanifolds of $M$) implies that the abelian group $
H_2(M)=H_2(M;\Zz )$ has no torsion. We shall view $H_2(M) $ as a
lattice in the $\Qq$-vector space $H_2(M;\Qq)=\Qq\otimes_{\Zz}
H_2(M) $. For any $x\in H_2(M;\Qq)$, there is an integer $n\geq 1$
such that $nx\in H_2(M)$.  Then $\Vert x\Vert_M= n^{-1} \min_\Sigma
\chi_- (\Sigma)\in \Qq $, where $\Sigma$ runs over all closed
oriented embedded surfaces in $M$ representing $
 nx  $. The number $\Vert x\Vert_M$ does not depend on the choice
of $n$ and is always realized by a certain $\Sigma$. Using surfaces
in $M$ with boundary on $\partial M$, one similarly defines the
Thurston   semi-norm   on $H_2(M,\partial M;\Qq)$.

\begin{lemma}\label{le2}    Let $j$ be  the
coefficient homomorphism $H_2(M;\Qq )\to H_2(M;\Qq/\Zz)$. Then
 $\ccc(j(x))\leq \Vert x\Vert_M$ for any $x\in H_2(M;\Qq )$.
\end{lemma}

\begin{proof} Let $\Sigma$ be a  closed oriented embedded
surface  in $M$ representing $  nx \in H_2(M )$ with $n\geq 3$. The
surface $\Sigma$ is an oriented folded surface with empty singular
set.  Consider  a folded surface $X=X(n)$ inside a   3-ball in $
M-\Sigma$, as constructed before Lemma \ref{le1}. The
  union $Z=X\cup \Sigma$ is a simple folded surface
representing $x$ and $i_Z=i_X=n$. By the definition of $\ccc $,
$$\ccc(j(x))\leq n^{-1}\, \chi_- (Z)=n^{-1}\,\chi_- ( \Sigma).$$
Therefore $\ccc(j(x))\leq \Vert x\Vert_M$.
\end{proof}

\begin{lemma}\label{le31}    Let $K$ be an oriented knot in $M$. Set  $N=M-K$ and
let  $\iota$ be  the  inclusion homomorphism $H_2(N;\Qq ) \to
H_2(M;\Qq )$. Let $j$ be the coefficient homomorphism $H_2(M;\Qq )
\to H_2(M;\Qq/\Zz)$.  Then for any simple folded surface $X$ in $M$
with $\sing (X)=K$ and any $y\in H_2(N;\Qq )$,
\begin{equation}\label{equ3}\ccc([ X ] +j \iota(y)) \leq
(i_X )^{-1} \chi_-(X)+ \Vert y\Vert_N.\end{equation}
\end{lemma}

\begin{proof} Set $n=i_X$ and let $k$
be a positive integer such that $ky\in H_2(N)\subset H_2(N;\Qq)$. It
is enough to prove that for any   closed oriented surface
$\Sigma\subset N$  representing  $nky$,
\begin{equation}\label{equ4}\ccc([ X ] +j \iota(y)) \leq n^{-1}
\chi_-(X)+ (nk)^{-1} \chi_-(\Sigma). \end{equation} This can be
reformulated in  terms of the simple folded surface $X^{(k)}$ as
$$\ccc([ X^{(k)} ] +j \iota(y)) \leq (nk)^{-1} (
\chi_-(X^{(k)})+   \chi_-(\Sigma)).$$ Therefore it is enough to
prove that for any simple folded surface $T$ in $M$  with $\sing
(T)=K$ and $i_T=nk$,
\begin{equation}\label{equ423}\ccc([ T ] +j \iota(y)) \leq (nk)^{-1}
(\chi_-(T)+   \chi_-(\Sigma)). \end{equation}

 Suppose first
that $T $ is {\it compressible} in $N=M-K$  in the sense that there
is an embedded closed 2-disk $D\subset N$ such that $T\cap
D=\partial D\subset T-K$ and the circle $\partial D $ does not bound
a 2-disk in $ T-K$. The surgery on $T$ along $D$ yields a simple
folded surface $T_D$ with $[ T_D] =[ T]$ and $\chi_-(T_D)<
\chi_-(T)$. Applying this procedure several times, we can reduce
\eqref{equ423} to the case  where $T$ is incompressible, i.e.,  $T$
admits  no disks $D$ as above. By the same reasoning, we can assume
that $\Sigma$ is incompressible in $N$ (it may be compressible in
$M$). The homology class $[ T ] +j \iota(y)\in H_2(M;\Qq/\Zz)$ is
represented by the 2-cycle $(nk)^{-1}  \,  T\cup  \Sigma  $ (mod
$\Zz$). Deforming   $\Sigma$ in $N$ so that it meets $T$
transversely and applying to $T\cup  \Sigma$ the usual cut and paste
technique, we can transform $T\cup  \Sigma$ into a simple folded
surface $Z$ with $\sing(Z)=\sing (T)=K$ and $i_Z=nk$. Clearly, $[ Z
]=[ T ]+j \iota( y)$. The folded surface $Z $ may have spherical
components (that is components homeomorphic to $S^2$) created from
pieces of $T-K$ and $\Sigma$ by cutting and pasting. One of these
pieces will necessarily be a 2-disk $D$ such that  either $D \subset
T-K$ and $D\cap \Sigma=\partial D$   or $D \subset \Sigma$ and
$D\cap (T-K)=\partial D$. In the first case the incompressibility of
$\Sigma$ implies that the circle $\partial D $ bounds a disk on
$\Sigma$. The surgery on $\Sigma$ along $D$ yields a surface
$\Sigma_+\approx \Sigma \amalg S^2$ homological to $\Sigma$ in $N$.
Then $\chi_-(\Sigma_+)=\chi_-(\Sigma)$ and the 1-manifold $T\cap
\Sigma_+$ has one component less than $T \cap \Sigma$. Similarly, if
$  D \subset \Sigma$, then the incompressibility of $T-K$ implies
that   $\partial D $ bounds a disk on $T-K$. The surgery on $T$
along $D$ yields a simple folded surface $T_+\approx T \amalg S^2$
such that  $[ T_+ ]=[ T ]$, $\chi_-(T_+)=\chi_-(T)$, and the
1-manifold $T_+\cap \Sigma$ has one component less than $T\cap
\Sigma$. Continuing in this way, we can reduce ourselves to the case
where $Z$   does not have spherical components  except the spherical
components of $T$   disjoint from $\Sigma$ and the spherical
components of $\Sigma$   disjoint from $T$. A similar argument
allows us to assume that
  the components of $Z-K$ are not disks except
the disk components of $T-K$   disjoint from $\Sigma$.
  Then the additivity  of the Euler characteristic under cutting and pasting implies   that
  $\chi_- (Z)  =
\chi_-(T)+\chi_- (\Sigma)$.  Therefore
$$\ccc([ T ] +j \iota(y)) \leq
 (nk)^{-1} \chi_-(Z) = (nk)^{-1} (\chi_-(T)+
 \chi_-(\Sigma)).$$
 This proves \eqref{equ423}, \eqref{equ4}, and \eqref{equ3}.
 \end{proof}

 \begin{thm}\label{le3}    If $M$ is compact, then there is a number $C >0$ (depending on $M$)  such that
$\ccc(x)\leq C$  for all $x\in H_2(M;\Qq/\Zz)$.
\end{thm}

\begin{proof} Set $d=d_M: H_2(M;\Qq/\Zz) \to H_1(M)$. Since the group
  $\Ima\, d=\Tors\, H_1(M) $ is
finite, it is enough to prove that for every   $u\in \Tors\,H_1(M)$,
the values of $\ccc $ on the elements of the set $d^{-1} (u)$ are
bounded from above.

Consider first the case $u=0$. Then $d^{-1}(u)= \Ima\, j$ where $j$
is   the coefficient homomorphism $ H_2(M;\Qq )\to H_2(M;\Qq/\Zz)$.
We need to prove that the values of $\ccc \circ j$ are bounded from
above. Since $M$ is compact, the group $H_2(M)$ is finitely
generated.
 Pick a basis $a_1,..., a_n $ in $H_2(M)$ and let $Q\subset
H_2(M;\Qq) $ be the cube consisting of the vectors $r_1 a_1+...
+r_na_n$ with
  rational non-negative $   r_1,...,r_n\leq 1$. The
    supremum $s=\sup_{x\in Q} \Vert x\Vert_M$  is a finite
number, because the Thurston semi-norm extends  to a continuous
semi-norm on $H_2(M;\Rr )$ and the closure of $Q$ in   $  H_2(M;\Rr
)$
  is   compact. We claim that $\ccc(j(x))\leq s$ for
any $x\in H_2(M;\Qq)$. Indeed, there is $a\in H_2(M) $ such that
$x+a\in Q$. Then   $j(x )=j(x+a)$ and  $\ccc(j(x))=\ccc(j(x+a))\leq
s$.

Consider now the case  $u\neq 0$. Pick an oriented knot $K\subset M$
representing $u$ and a simple folded surface $X$ in $M$ with $\sing
(X)=K$.  Then $d^{-1}(u)=\{[ X ] +j \iota(y)\}_y$ where  $\iota$ is
the inclusion homomorphism $H_2(M-K;\Qq ) \to H_2(M;\Qq )$ and $y$
runs over $H_2(M-K;\Qq)$. The rest of the argument goes as  in the
case $u=0$ using Lemma \ref{le31}.\end{proof}

\subsection{Semi-continuity} For compact $M$, the  group
$H_2(M; \Qq/\Zz)$ has a natural topology    as follows. The image of
the coefficient  homomorphism $j: H_2(M;\Qq)\to H_2(M; \Qq/\Zz)$ can
be identified with the quotient   $H_2(M;\Qq)/H_2(M)$. Provide
$\Ima\, (j)$ with the quotient topology induced by the standard
topology in the finite dimensional $\Qq$-vector space $H_2(M;\Qq)$.
This extends to a topology in $H_2(M; \Qq/\Zz)$ by declaring   a set
$U\subset H_2(M; \Qq/\Zz)$ open if $(a+U)\cap \Ima\, (j)$ is open in
$\Ima\, (j)$ for all $a\in H_2(M; \Qq/\Zz)$. Recall that an
$\Rr$-valued function $f$ on a topological space $A$ is {\it upper
semi-continuous} if for any point $a\in A$ and any real
$\varepsilon>0$, there is is a neighborhood $U\subset A$ of $a$ such
that $f(U)\subset (-\infty, f(a)+ \varepsilon)$.

\begin{lemma}\label{le4}    For compact $M$, the function $\ccc=\ccc_M $ is upper semi-continuous.
\end{lemma}

\begin{proof} Let $a\in H_2(M; \Qq/\Zz)$ and $\varepsilon>0$. Let
$X$ be a simple folded surface in $M$ representing $a$ and such that
$(i_X)^{-1} \chi_-(X)\leq \ccc(a ) +\varepsilon/2$. Set $K=\sing
(X)$ and $N=M-K$.  Let   $\iota:H_2(N;\Qq )\to H_2(M;\Qq ) $ be the
inclusion homomorphism. Put
$$V=\{y\in H_2(N;\Qq )\, \vert \,  \Vert y\Vert_N <
\varepsilon/2\}.$$ The set $V$ is open in $H_2(N;\Qq )$ since the
Thurston  norm is continuous. The set $\iota(V)$ is open in $
H_2(M;\Qq )$ since $\iota$ is an epimorphism. The set $  j \iota(V)$
is open in $\Ima\, (j)$ by  definition of the topology in $\Ima\,
(j)$. Finally, the set $U=a+j \iota(V)$ is an open neighborhood of
$a$ in $H_2(M; \Qq/\Zz)$ by   definition of the topology in $H_2(M;
\Qq/\Zz)$. By \eqref{equ3},  $\ccc(U) \subset (-\infty, \ccc(a)+
\varepsilon)$. Hence $\ccc $ is upper semi-continuous.
\end{proof}

\section{Estimates  from below: the case $b_1  \geq 1$}\label{section:3k}

In this section we give an estimate from below for the functions
$\ccc =\ccc_M$ and $\CCC=\CCC_M $ of a 3-manifold $M$  with non-zero
first Betti number  $b_1(M)$.

We begin with preliminaries on group rings and abelian torsions of
3-manifolds.

\subsection{Preliminaries}\label{tor}  Let $H$ be a finitely generated abelian
group written in  multiplicative notation.  Any element $a$ of the
group ring $\Qq [H]$ expands uniquely in the form
   $a=\sum_{h\in H} a_h h  $, where $a_h\in \Qq$ and $a_h=
0$ for all but finitely many $h$. We say that an element $h\in H$ is
{\it $a$-basic} if $a_h\neq 0$. The (finite) set of $a$-basic
elements of $H$ is denoted  $B_a$. The element $\sum_{h\in \Tors\,
H} h$ of $\Qq [H]$ will be denoted $\Sigma_H$. Clearly,
$B_{\Sigma_H}=\Tors\, H$.

The classical ring of quotients  of   $\Qq [H]$ that is, the
(commutative) ring obtained by inverting all non-zero-divisors of
$\Qq [H]$ is denoted $Q(H)$. It  is known that   $\Qq[H]$ splits as
a direct sum of domains. Therefore $Q(H)$ splits as a direct sum of
fields and the natural ring homomorphism  $  \Qq[H] \to Q(H)$ is an
embedding. We identify   $\Qq[H]$ with its image under this
embedding. Note that if $H$ is a finite abelian group, then
$Q(H)=\Qq[H]$.

 Let $M$ be a compact connected  3-manifold.   From now on,
  we use multiplicative notation for the group operation in
$H=H_1(M)$. In particular, the neutral element of $H$ is denoted
$1$. The manifold $M$ gives rise to a {\it maximal abelian torsion}
$\tau(M)$ which is an element of   $Q(H) $ defined up to
multiplication by $-1$ and elements of $H$,   see \cite{tu, ni}. If
$b_1(M)\geq 2$, then all representatives of $\tau(M)$ belong to $\Zz
[H]\subset \Qq[H] \subset Q(H)$. We express this by writing $\tau(M)
\in \Zz [H]$. If $b_1(M)=1$ and $\partial M\neq \emptyset$, then
$\tau(M)\in \Zz [H]+\Sigma_{ H} \cdot Q(H)$. This implies that
$(h-1)\,\tau(M)\in \Zz [H]  $ for all $h\in \Tors\, H$ (indeed
$(h-1)\Sigma_{ H}=0$).

  If $M$ is oriented and
  $b_1(M)\geq 2$, then  the Thurston  semi-norm $\Vert
\cdot \Vert_M$ on $H_2(M,\partial M;\Qq)$ can be estimated in terms
of   $\tau(M)$ as follows (see \cite{tu}):  for any $s\in H_2
(M,\partial M;\Qq)$ and any representative $a\in \Zz [H]$ of
$\tau(M)$,
\begin{equation}\label{TT} \Vert s\Vert_M \geq \max_{h,h'\in B_a} \vert   h \cdot
s -   h'\cdot s  \vert, \end{equation} where    $h\cdot s \in \Zz$
is the intersection index of $h$ and $s$. Note that the right hand
side of  \eqref{TT} does not depend on the choice of  $a$ in
$\tau(M)$.

 \subsection{An estimate for $\ccc_M$} The function $ \ccc $ will be  estimated
   in terms of spans of   subsets of $\Qq/\Zz$.
  The  {\it span} $\spn(A)$ of a finite  set $A\subset \Qq/\Zz$  is a rational  number
  defined as the minimal length of an interval in $\Qq/\Zz$
containing
  $A$, that  is
the minimal   rational number $t\geq 0$ such that for some $r\in \Qq
$, the projection of the set  $[r,r+t]\cap \Qq$ into $\Qq/\Zz$
contains $A$. Clearly, $1>\spn(A) \geq 0$ and $\spn(A)=0$ if and
only if $A$ is empty or has only one element.

Given an oriented 3-manifold $M$ and a homology class $x\in
H_2(M;\Qq/\Zz)$, we set for any   $a\in \Qq [H_1(M)]$,
 $$\spn_x(a) =\spn (\{ h\cdot x\}_{h\in B_a}),$$
 where    $h\cdot x \in \Qq/\Zz$
is the intersection index of $h$ and $x$. Clearly, $1>\spn_x(a) \geq
0$.

 \begin{thm}\label{esb1}  Let $M$  be a
compact connected oriented 3-manifold     with $b_1(M)\geq 1$. Set
$H=H_1(M)$ and let $\tau \in Q(H)$ be a representative of the
torsion $\tau(M)$. Let $x\in H_2 (M;\Qq/\Zz)$ and $u=d_M(x)\in H$.
Then $(u-1)\, \tau \in \Zz [H]$ and
\begin{equation}\label{ineq}  \ccc_M(x) \geq \spn_x( (u-1)\, \tau).
\end{equation}  \end{thm}

\begin{proof} If $b_1(M)\geq 2$, then $\tau \in \Zz [H]$ and
  $(u-1)\, \tau \in \Zz [H]$.  The inclusion  $u\in
\Tors\, H$  and the remarks in Section \ref{tor} imply that $(u-1)\,
\tau \in \Zz [H]$ for $b_1(M)=1$ as well.

We   prove \eqref{ineq}. Let $X $ be a simple folded surface in $M$
representing $x$. The knot $\sing (X) \subset M$ endowed with
orientation induced from the one on $X$ represents  $u\in \Tors\, H
$. Let $E$ be the exterior of this knot in $M$.   The homological
sequence of the pair $(M,E)$ and the inclusion $u\in \Tors\, H $
imply that $b_1(E)\geq b_1(M)+1\geq 2$. Therefore $\tau(E)\in
\Zz[H_1(E)]$. Pick a representative $a\in \Zz[H_1(E)]$ of $\tau(E)$.
Denote by $\iota$ the inclusion homomorphism $H_1(E)\to H_1(M)=H$
and denote $\iota_\ast$ the induced ring homomorphism $\Zz
 [H_1(E)] \to \Zz [ H  ]$. By \cite{tu}, Theorem VII.1.4,  $\iota_\ast(a)=
  (u-1)\,b$ where $b$ is a representative of $\tau(M)$.
 Note that the right hand side of \eqref{ineq} does not depend on the choice of $\tau$
  in $\tau(M)$. Therefore without loss of generality we can assume that
 $ \tau=b$.

Deforming, if necessary,  $X$ in $M$, we can assume  that $S=X\cap
E$ is the complement in $X$ of a regular neighborhood of $\sing
(X)$. Then $S $ is a proper surface in $E$ and $\chi_- (X)  = \chi_-
(S)$. The orientation of $\Int (X)$ induces an orientation of $S$.
The oriented surface  $S $ represents a relative homology class
$s\in H_2(E,
\partial E)$. By \eqref{TT},  $$ \chi_- (X)  =\chi_- (S) \geq \max_{h,h'\in B_a} \vert
  h \cdot s
-   h'\cdot s  \vert   $$ where $B_a \subset H_1(E) $ is the set of
$a$-basic elements.  Let $r\in \Qq$ be the minimal element of the
set $\{h\cdot s\}_{h\in B_a}$. Then
$$\{h\cdot s\}_{h\in B_a} \subset [r, r+\chi_- (X)].$$
Denote the projection $\Qq\to \Qq/\Zz$ by $\pi$. Observe that  for
any $h\in H_1(E)$, $$\iota (h)\cdot x = \pi \left (\frac{h\cdot
s}{i_X} \right ).$$   Therefore
$$\{\iota(h)\cdot x\}_{h\in B_a} \subset \pi \left ( \left [\frac{r}{i_X},\, \frac{r+\chi_- (X)}{i_X} \right ]\right ).$$
 The equality $\iota_\ast(a)=
  (u-1)\,\tau $ implies that
  $B_{(u-1)\,\tau}  \subset \iota  (B_a) $. Hence
$$\{g\cdot x\}_{g\in B_{(u-1)\,\tau}} \subset \{\iota (h)\cdot x\}_{h\in B_a} \subset
 \pi \left ( \left [\frac{r}{i_X},\, \frac{r+\chi_- (X)}{i_X} \right ]\right ).$$
  Therefore
$$   \spn_x( (u-1)\, \tau) \leq  (i_X)^{-1}\,  \chi_- (X) .$$
Since this holds for all simple folded surfaces $X$ representing
$x$, we have \eqref{ineq}.
\end{proof}

\subsection{An estimate for $ \CCC_M$} Let $M$ and $H$ be as in Theorem \ref{esb1}.
To estimate the function $\CCC_M:\Tors \, H \to \Qq/\Zz$, we need
the linking form $L_M:\Tors\, H\times \Tors\, H\to \Qq/\Zz$ of $M$.
It is    defined by $L_M(h,g)=  h \cdot x\in \Qq/\Zz$ where $x$ is
an arbitrary element of $H_2(M;\Qq/\Zz)$ mapped to $g$ by the
boundary homomorphism $ d:H_2(M;\Qq/\Zz) \to  H$. The pairing $L_M$
is well defined, bilinear, and symmetric.

Given     $u\in \Tors\, H $ and $a\in \Qq [H ]$,   set
 $$\spn_u(a) =\spn (\{ L_M(h,u) \}_{ h\in B_a\cap \Tors\, H }).$$
 Clearly, $\spn_x(a) \geq \spn_{d(x)} (a)$ for any $x\in
 H_2(M;\Qq/\Zz)$ and any $a\in \Qq [H ]$. This   and
Theorem \ref{esb1} imply that,   under the conditions of this
theorem,
\begin{equation}\label{labCCC} \CCC_M(u)\geq \spn_u( (u-1)\,
\tau),\end{equation} for  any $u\in \Tors\, H $ and any
representative $\tau$ of $\tau(M)$.  Generally speaking,  the
right-hand side of \eqref{labCCC} depends on the choice of $\tau$.

\subsection{Remark} Estimate \eqref{TT} strengthens the McMullen
\cite{mc} estimate of the Thurston  norm via the Alexander
polynomial. For recent more general estimates of this type, see
Friedl \cite{fr}.

\section{Estimates  from below: the case of
$\Qq$-homology spheres}\label{section:4k}

For $\Qq$-homology spheres, the functions $\ccc$ and $\CCC$ contain
the same information and it is enough to give an estimate for
$\CCC$. We begin with preliminaries on  refined  torsions  and
$\Qq$-homology spheres, referring for details to \cite{tu}, Chapters
I and   X.

\subsection{Refined torsions} The maximal abelian torsion $\tau(M)$
of a compact connected 3-manifold $M$ admits a refinement
$\tau(M,e,\omega)\in Q(H_1(M))$ depending on an orientation $\omega$
in the vector space $H_\ast(M;\Qq)=\oplus_{i\geq 0} H_i (M;\Qq)$ and
 an Euler structure $e$ on $M$. An Euler structure on $M$ is
determined by a non-singular vector field on $M$ directed outside on
$\partial M$. Two such vector fields determine the same Euler
structure if for a point $x\in \Int (M)$, the restrictions of these
fields to $M-\{x\}$ are homotopic in the class of non-singular
vector field on $M-\{x\}$ directed outside on $\partial M$. The set
of Euler structures on $M$ is denoted $\Eul (M)$. This set admits a
canonical free transitive action of the group $  H_1(M)$. The
torsion $\tau(M,e,\omega)$ satisfies $\tau(M,he,\pm \omega)= \pm h\,
\tau(M,e,\omega)$ for any $e\in \Eul (M), h\in H_1(M) $. The
unrefined torsion $\tau(M)$ is just the set $\{\pm
\tau(M,e,\omega)\}_{e\in \Eul (M)}$. If $\partial M=\emptyset$, then
the set $\Eul (M)$ can be identified with the set of
$Spin^c$-structures on $M$.

\subsection{Homology spheres}    Let $M$ be an oriented
 3-dimensional $\Qq$-homology sphere. Denote
 $\omega_M$ the orientation in $H_\ast(M;\Qq)= H_0 (M;\Qq)\oplus H_3
(M;\Qq)$ determined by the following basis: (the homology class of a
point,   the fundamental class of $M$).

The   group $H=H_1(M)$ is finite and the linking form $L_M:H\times
H\to \Qq/\Zz$ is non-degenerate in the sense that the adjoint
homomorphism $H\to \Hom (H, \Qq/\Zz)$ is an isomorphism. Recall that
we use multiplicative notation for the group operation in $H$. Every
Euler structure $e\in \Eul(M)$ determines a torsion $\tau(M,e,
\omega_M)\in Q(H)=\Qq[H]$. The linking form $L_M$ can be computed
from this torsion by
\begin{equation}\label{stra3t}L_M(h,g)=
   - \pi(((1-h)(1-g)\,\tau(M,e, \omega_M))_1) \in
\Qq/\Zz \end{equation}
    for all $h,g \in H$, where  $\pi$  is the projection $\Qq\to
\Qq/\Zz$ and for any $ a\in \Qq[H]$, the symbol $a_1\in \Qq$ denotes
the coefficient of the neutral element   $1\in H$ in the expansion
of $a$ as a formal linear combination of elements of $H$ with
rational coefficients. The Euler structure $e$ determines a function
$q_e:H\to \Qq/\Zz$ by
\begin{equation}\label{stra}q_e(u)=\pi(((1-u)\,\tau(M,e, \omega_M))_1),\end{equation} for any $u\in H$. It follows from
\eqref{stra3t}, \eqref{stra} that $q_e$ is
  quadratic in the sense that
$q_e(hg)=q_e(h)+q_e(g)+L_M(h,g)$ for all $h,g\in H$. Formula
\eqref{stra} also implies that
\begin{equation}\label{stra3hj}q_{he} (u)=q_e(u)+ L_M(h,u), \end{equation}
for any $h\in H$.

 If $u\in H$ has order $n$
(i.e., $n$ is the minimal positive integer
  such that $u^n=1$), then by \cite{tu}, Section X.4.3 there is a unique residue $K(e,u) \in
\Zz/2n \Zz$ such that
\begin{equation}\label{resid}q_e(u)= \frac{K(e,u)}{2n} +\frac{1}{2}\, ({\text {mod}}\, \Zz).\end{equation}
 Formula \eqref{stra3hj} implies that the residue  $K(e,u)\,
({\text {mod}}\,2)$ does not depend on $e$. We say that $u$ is {\it
even} if this residue is 0 and {\it odd} if it is 1.

Every homology class
  $u\in H$  gives rise to a group  $$ G =G_u= \{g\in H\,\vert \,
L_M(u,g)=0\}\subset H.$$ The non-degeneracy of $L_M$ implies that
the quotient $ H/G $  is a finite cyclic group whose order  is equal
to the order, $n$, of $u$ in $H$. Moreover, there is an element
$v=v_u \in H$ such that $L_M(u,v)= n^{-1}\, ({\text {mod}}\, \Zz)$.
Such $v$ is determined by $u$ uniquely up to multiplication by
elements of $G$. The inclusion $v^n\in G$ implies that the order,
$p$,  of $v$ is divisible by $n$ (in particular, $p\geq n$). Set
\begin{equation}\label{iop} \alpha_v\,=\, \frac {1+2v+3v^2+...+pv^{p-1}}{p}\, - \, \frac
{p+1}{2} \cdot \frac {1+v+...+v^{p-1}}{p} .
\end{equation}
This element of $\Qq[H]$   can be uniquely characterized by the
following property: for any ring homomorphism $\varphi$ from
$\Qq[H]$ to a field, $\varphi(v)=1 \Rightarrow \varphi(\alpha_v) =0$
and $\varphi(v)\neq 1 \Rightarrow \varphi(\alpha_v) =(\varphi
(v)-1)^{-1}$. In \cite{tu}  we used the notation $(v-1)^{-1}_{\text
{par}}$ for $\alpha_v$.

 \begin{thm}\label{esb3}  Let $M$  be an
oriented 3-dimensional $\Qq$-homology sphere.  Let   $u$ be an
element of $ H=H_1(M)$ of order  $n\geq 1$. Set $G =\{g\in H\,
\vert\, L_M(u,g)=0\}$ and  $\Sigma_G=\sum_{g\in G} g\in
\Zz[G]\subset \Zz[H]$. Pick any $v\in H$ such that $L_M(u,v)=n^{-1}
\,({\text {mod}}\,  1)$. For  $e \in \Eul (M)$, set
$$a_e(u)=(u -1)\,  \tau
 (M,e, \omega_M
) -  \,  \frac{  v^{K(e ,u) /2} (v+1)}{2} \,  \alpha_v\, \Sigma_{G}
\in \Qq[H],$$ if
 $u$ is even  and
 $$a_e(u)=(u -1)\,  \tau
 (M,e, \omega_M
) - v^{(K(e ,u)+1) /2} \, \alpha_v\,   \Sigma_{G}\in \Qq[H],
$$  if $u$ is odd. Then for any $e
\in \Eul (M)$,
 $$\CCC_M (u)  \geq \spn_u   (a_e(u))= \spn (\{ L_M(h,u) \}_{ h\in B_{a_e(u)} })  .$$ \end{thm}

\begin{proof}     If $u$ is even (resp.\ odd), then $K(e,u)\in \Zz_{2n}$ is even
(resp.\ odd). Therefore the power of $v$ in the definition of
$a_e(u)$ is well defined up to multiplication by $v^n$. However,
$v^n\in G$ and    $v^n\Sigma_G=\Sigma_G$. Therefore the right hand
sides of the   formulas for $a_e(u)$ are well defined. If $v'$ is
another element of $H$ such that $L_M(u,v')= n^{-1} \,({\text {mod
}} 1)$, then $v' \in v\,G$ and $v^k\Sigma_G= (v')^k \Sigma_G$ for
all $k\in \Zz$. Therefore $a_e(u)$ does not depend on the choice of
$v$. It is easy to see that $a_{he}(u)= h \,  a_e(u)$ for all $h\in
H$. Therefore the number $\spn_u (a_e(u))$ does not depend on $e$.

 Let $X\subset M$ be a simple folded surface
representing the 2-homology class $x=d_M^{-1}(u)\in H_2
(M;\Qq/\Zz)$. The knot $K=\sing (X)  $ with orientation induced from
the one on $X$ represents  $u\in H_1(M) $. Let $E$ be the exterior
of $K$ in $M$.
  Clearly $b_1(E)=1$. Fix an orientation $\omega$ in $H_\ast (E;\Qq)$  and
   an   Euler structure $e_K$ on $E$. The torsion $\tau (E,e_K,\omega)\in Q(H_1(E))$
   can be canonically expanded as a sum of a certain
 $[\tau]=[\tau]  (E,e_K, \omega) \in  \Qq[H_1(E)]$  with an element
 of $Q(H_1(E))$ given by an
explicit formula using solely   $\omega$  and the Chern class
 of $e_K$, see \cite{tu}, Section II.4.5. The
 inclusion homomorphism $\Qq[H_1(E)]\to \Qq[H_1(M)]$ sends $[\tau]$  to
 $\pm a_e(u)$ for some $e\in \Eul (M)$, see \cite{tu}, Formula X.4.d.
 The inequality \eqref{TT} holds for any $s\in H_2(E,\partial E;\Qq)$ and
$a=[\tau]  $, see \cite{tu}, Chapter IV. The rest of the argument
goes
 as the proof of Theorem \ref{esb1}
 with $\tau$ replaced by $[\tau]  $. This gives
 $ (i_X)^{-1}  \chi_- (X)  \geq \spn_x( a_e(u)) =\spn_u( a_e(u))$.
Since this holds for all $X$ representing $x$, we have
$\CCC_M(u)=\ccc_M(x)   \geq \spn_u   (a_e(u))$.
\end{proof}

\subsection{Remarks} 1. Let $ \frac{1}{2}\Zz$
be the additive group of integers and half-integers. In  Theorem
\ref{esb3}, $a_e(u)\in \Zz[H]$ if $u$ is even and $a_e(u) \in
\frac{1}{2}\Zz [H ]$ if $u$ is odd.  This follows from the proof of
this theorem and the inclusion  $[\tau]   \in \Zz [H_1(E)]$ if $u$
is even and $[\tau]   \in \frac{1}{2}\Zz [H_1(E)]$ if $u$ is odd.

2. It is  proven in \cite{dm} that the   function $q_e:H\to \Qq/\Zz$
derived from the torsion   coincides with the quadratic function
defined geometrically in   \cite{lw}, \cite{de}.

\section{Examples}\label{section:36k}

\subsection{Lens spaces} The computation of the abelian torsions
 for  the lens space $M=L(p,q)$   goes back to K. Reidemeister, see, for instance,
\cite{tu2} for an introduction to the theory of torsions. Let $t$,
$t^q$ be the generators of $H=H_1(M)$ represented by the core
circles of the two solid tori forming $M$. For an appropriate choice
of an orientation on $M$ and an
 Euler structure $e$ on $M$,   we
have
 $\tau(M,e, \omega_M)=\alpha_t\, \alpha_{t^q}$,  where $\alpha_v\in
\Qq[H]$ is defined  by  \eqref{iop} for any $v\in H$.   This  allows
us to compute $a_e(u)$ for any $u\in H$ and to apply Theorem
\ref{esb3}.  We give here
  a few    examples.

Consider   the lens space $M=L(5,1)$.  By Sections \ref{der}, $\CCC
(t^{4})=\CCC (t)=\CCC (1)=0$ and $\CCC(t^2)=\CCC(t^3)$. We show that
$\CCC (t^2)\geq 1/5$. We have
$$\alpha_t=\frac{-2-t+t^3+2t^4}{5}.$$
Then
$$\tau=\tau(M,e, \omega_M)=\alpha_t^2=\frac{t+t^2-2t^4}{5}.$$
A direct computation shows that $L_M(t,t)= (-(1-t)^2 \tau)_1=1/5$
and $q_e(t^2)=((1-t)\,\tau)_1=0$. Note that $u=t^2$ has order $5$ in
$H$. From \eqref{resid},
  $K(e,u)=5\, ({\text {mod}}\, 10)$. Therefore $u$ is odd. The
 associated group $G_u$ is trivial, $v=v_u=t^3$, and
$$a_e(u)= (u -1)\,  \tau
   - v^3 \alpha_v=t^4-t.$$
   Since
    $L_M(t^4,u)=3/5 \, ({\text {mod}}\, 1 )$ and $L_M(t,u)=2/5 \, ({\text {mod}}\, 1 )$,
    the span of the set
   $\{ L_M(h,u) \}_{ h\in B_{a_e(u)} }$ is equal to $1/5$. By
   Theorem \ref{esb3}, $\CCC
(t^2)\geq 1/5$. In this example, the   function $\CCC:H\to \Rr_+ $
takes non-zero values only on $t^2$ and $t^3$. This function is
 non-homogeneous and does  not  satisfy  the triangle inequality.

Consider   the lens space $M=L(6,1)$. Then
$$ \alpha_t =\frac{-5-3t-t^2+ t^3+3t^4+5t^5}{12}, $$
$$\tau=\alpha_t^2=\frac{-5+13t+19t^2+13t^3-5t^4-35t^5}{72}, $$
and $L_M(t,t)=1/6$. For  $u=t^2$, the computations similar to the
ones above give $q_e(u)=0 \, ({\text {mod}}\, 1 )$, $ K(e,u)=3
({\text {mod}}\, 6)$, $G_u=\{1, t^3\}$, $v_u=t$, and $a_e(u)=t^5-t$.
   Theorem \ref{esb3} yields  $\CCC
(t^2)\geq 1/3$. For  $u=t^3$, we similarly obtain   $q_e(u)=3/4 \,
({\text {mod}}\, 1 )$, $K(e,u)=1\, ({\text {mod}}\, 4)$, $G_u=\{1,
t^2, t^4\}$, $v_u=t$, and  $a_e(u)=t^5-t^2$.
   Theorem \ref{esb3} yields  $\CCC
(t^3)\geq 1/2$.

\subsection {Surgeries on knots}\label{olp}  Let   $L$ be an oriented knot in
an oriented 3-dimensional $\Zz$-homology sphere $N$.
  Let   $M$ be the closed oriented
  3-manifold obtained by surgery on $N$ along  $L $ with
framing $  p \geq 2$. Let $u \in H=H_1(M) $ be the homology class of
the meridian of $L$ whose linking number with $L$ is   $+1$.
Clearly, $H$ is a cyclic group of order $p$ with generator $u$ and
$L_M(u, u)=   p^{-1} \,({\text {mod}}\, 1)$. We
 explain now how to estimate  $\CCC (u)$ in terms of the Alexander polynomial of
$L$. We will see that in  some cases this estimate is exact.

Recall   that the {\it span} $\spn(\Delta)$ of a non-zero Laurent
polynomial $\Delta=\sum_i a_i t^i \in \Zz [t^{\pm 1}]$ is the number
$ {\max \{i\,\vert\, a_i\neq 0\}}- {\min(\{i\,\vert\, a_i\neq
0\})}$. Let $\Delta=\Delta_{L }(t)$ be the Alexander polynomial of
$L $ normalized so that $\Delta(t^{-1})=\Delta(t)$ and $\Delta(1)=
1$. Expand $\Delta(t)=1+ (t-1) \, \beta(t)$ where $\beta(t)\in
\Zz[t^{\pm 1}]$. We claim  that the expression $a_e(u)\in \Qq [H]$
defined in
 Theorem \ref{esb3} is equal to $\beta(u)$ for an appropriate Euler
structure $e$ on $M$. By Theorem \ref{esb3}, this will imply that
$\CCC(u)\geq \spn_u (\beta(u))$. For example, if $p \geq 2 \,\spn
(\beta)$, then $\spn_u (\beta(u)) = p^{-1} {\spn (\beta)} = p^{-1}
({\spn (\Delta) -1})  $. Therefore $\CCC(u)\geq p^{-1} ({\spn
(\Delta) -1}) $. On the other hand,
  by Section \ref{coae}.3,  $\CCC(u)\leq p^{-1} (2g-1)$,
where $g$ is the genus of $K$. In particular, if $\spn
(\Delta)=2g>0$ (for instance, if $K$ is a non-trivial fibred knot)
and $p \geq 4g-2$, then $\CCC(u)= p^{-1} (2g-1)$.

We now verify the claim above. Set $\tau=\alpha_u^2\, \Delta(u) \in
\Qq[H]$. It is easy to deduce from the multiplicativity of the
torsions that $\tau(M,e, \omega_M)=\tau$ for a certain orientation
on $M$ and a certain Euler structure $e$ on $M$ (for details, see
\cite{tu}, Formula X.5.e).
 Set $\sigma=1+u+u^2+\cdots + u^{p -1}\in \Zz [H]$.
 Clearly,  $\sigma u^k =\sigma   $ for any  integer $k$. Therefore
 for any integer 1-variable polynomial $f$, the product $\sigma
 f(u)$ is equal to $ \aug(f)\, \sigma$ where $\aug (f)=f(1)$ is the sum of
 coefficients of $f$.
 Since   $\aug(\alpha_u)=0$, we have $\sigma \alpha_u=0$. A direct computation shows  that $( 1-u)
\alpha_u= \sigma/p-1$. Hence
 $$ (1-u)\,
\tau =(1-u)\, \alpha^2_u \,\Delta(u)= ( \sigma/p-1)\,\alpha_u
\,\Delta(u)=-\alpha_u \,\Delta(u) $$
$$= -\alpha_u +\alpha_u  \,(1-u)\, \beta(u)=
-\alpha_u+(\sigma/p-1) \,\beta(u)=-\alpha_u -\beta(u),$$ where we
use the equality $\aug(\beta )=0$ which follows from the symmetry of
$\Delta$. Thus,
$$q_e(u)=((1-u) \tau)_1=-(\alpha_u)_1= (p-1)/2p \,({\text {mod}}\, 1).$$  Formula
\eqref{resid} implies that $K(e,u)=-1\,({\text {mod}}\, 2p)$. In
particular, $u$ is odd.

We also have
 $$(1-u)^2 \tau=(1-u)(-\alpha_u-\beta(u))= 1-\sigma/p -(1-u) \beta(u)  .$$
Hence $$L_M(u,  u)= -((1-u)^2 \tau)_1 = p^{-1}\, ({\text {mod}}\,
1).$$ This shows that the orientation of $M$ chosen so that
$\tau(M,e, \omega_M)=\tau$ is actually the one induced   from the
orientation on $N$. The equality $L_M(u,  u)= p^{-1}\, ({\text
{mod}}\, 1)$ implies that $v_u=u$ and $G_u=1$. We conclude that
$$a_e(u)=(u-1)\,\tau-\alpha_u=\alpha_u +\beta(u)
 -\alpha_u=  \beta(u).$$

\subsection{Surgeries on 2-component links} Let $M$ be a closed
oriented 3-manifold obtained by surgery on a 2-component oriented
link $L=L_1\cup L_2 $ in an oriented 3-dimensional $\Zz$-homology
sphere $N$. Suppose that   the linking number of $L_1, L_2$ in $N$
is $0$, the framing  of $L_1$ is $p\neq 0$, and the framing   of
$L_2$ is $0$. Then $H=H_1(M)=(\Zz/p\Zz)u_1 \oplus \Zz u_2$, where
   $u_i\in H$ is the homology class of the meridian of
$L_i$
 whose linking number with $L_i$ is
$+1$, for $i=1,2$. The Alexander polynomial of   $L$ has the form
$$\Delta_L (t_1, t_2)= f(t_1,t_2) (t_1-1) (t_2-1)$$
for some Laurent polynomial  $f(t_1,t_2)\in \Zz [t_1^{\pm 1},
t_2^{\pm 1}]$. Both $\Delta_L$ and $f$ are defined only up to
multiplication by $-1$ and monomials on $t_1, t_2$. By  \cite{tu},
Formula VIII.4.e, the torsion $\tau(M)$ is represented by
$$\tau=f(u_1,u_2) \pm  \Delta_{L_2}(u_2)\, u_2^n \,(u_2-1)^{-2} \Sigma_H \in Q(H)  $$
for an appropriate sign  $\pm$ and an integer  $n $, both depending
on the choice of $f$. Here $\Delta_{L_2}$ is the Alexander
polynomial of   $L_2$ normalized as in Section \ref{olp}. Pick $x\in
H_2 (M;\Qq/\Zz)$ and set $u=d(x)\in \Tors\, H$. Since
$(u-1)\Sigma_H=0$, Theorem \ref{esb1} implies that
\begin{equation}\label{ineq2}  \ccc (x) \geq \spn_x( (u-1)\, f(u_1,u_2)).
\end{equation}
For sufficiently big $p$,  the span on the right hand side does not
depend on $p$.

Note   another curious phenomenon. Suppose for simplicity that
$f(t_1, t_2)=1$ (a constant polynomial). Then $\ccc(x) \geq \spn_x(
u-1)$.
 If
$u=d(x)\neq 1$, then the set $B_{u-1 }\subset H$ consists of two
elements $u , 1$ and  $$\spn_x(u-1) =\spn (\{ u\cdot x, 0\})=\spn
(\{ L_M(u,u), 0\}).$$ For $u=u_1^k$ with $k\in \{0, 1,..., n-1\}$,
we have  $L_M(u,u) =k^2/n \,({\text {mod}} \,1)$. For
$k<\sqrt{n/2}$, we obtain $\spn (\{ L_M(u,u), 0\})= k^2/n $. Thus
$\CCC(u_1^k) \geq k^2/n $.  This  suggests that   the number
$\CCC(u_1^k)$, considered as a function of $k$, may behave like
  a quadratic function for small values of $k$.

\section{Miscellaneous}\label{section:3hjk6k}

\subsection{Quasi-simple folded surfaces} One can use a larger class of folded surfaces to
represent 2-homology classes. Let us call a folded surface $X$ {\it
quasi-simple} if it is oriented, $\sing(X)\neq \emptyset$,  and the
indices of all components of $ \sing (X)$ in $X$ are equal to each
other and non-zero. Denote the common value of these indices $i_X$.
  In particular, simple folded surfaces are
quasi-simple.

For a  quasi-simple folded surface $X$ in a 3-manifold $M$,
  the   2-chain
  $ (i_X)^{-1} X$
  is a 2-cycle mod  $\Zz$
   representing   a homology class  $[ X ] \in H_2(M;\Qq/\Zz)$. We claim that
   \begin{equation}\label{eq245} \ccc([X])\leq i_X^{-1}
 \,\chi_- (X) +b_0(\sing(X)) -1 \end{equation}  where $b_0(\sing(X))$ is the number
of components of $\sing(X)$. Indeed,   $X$   can be modified in a
neighborhood of   $\sing (X)$ so that each point of $\sing(X)$ is
adjacent   to exactly $i_X$ local branches of $\Int (X)$ (which then
induce the same orientation on $\sing(X)$). Let $\Gamma$ be a graph
with two vertices and $  i_X  $ edges connecting these vertices.
Given an embedded arc in $M$ with endpoints on different components
of $\sing (X)$ and with interior in   $M-X$, we can modify $X$ by
cutting it out along $\sing (X)$ near the endpoints and gluing in
$\Gamma\times [0,1]$  along the arc. This gives a quasi-simple
folded surface, $Z$, such that
$$b_0(\sing(Z))=b_0(\sing(X))-1,\,  i_Z=i_X,\,  [Z]=[X],\, {\text { and}}\,  \chi_-
(Z)\leq \chi_-(X)+i_X.$$  Modifying $X$ in this way, we can reduce
ourselves to the case where $\sing(X)$ is connected. In this case
\eqref{eq245} follows from the definition of $\ccc$. It may happen
that  there are no distinct components of $\sing(X)$  connected by
an arc with interior in $M-X$. This occurs if each   arc joining
distinct components of $\sing(X)$ has to cross the closed 2-manifold
$X_0 $ formed by  the components of $X$ disjoint from $\sing (X)$.
To circumvent this obstruction, we first modify
  $X_0$ so that $X-X_0$ is contained in a  connected component of
  $M-X_0$, cf.\    \cite{tu}, p.\ 60.

Formula \eqref{eq245} implies that for any $x\in H_2(M;\Qq/\Zz)$,
   \begin{equation}\label{eq1+9} \ccc(x)=\inf_X \,  \left (\frac
{\chi_- (X) }{i_X} +b_0(\sing(X))\right ) -1,\end{equation} where
$X$ runs over all quasi-simple folded surfaces
 in $M$    representing~$x$.

\subsection{Coverings}   Let $M$ be a compact   oriented 3-manifold and
 $p:\widetilde M\to M$ be an $n$-fold (unramified)
 covering. Let $p^*:H_2(M;\Qq/\Zz) \to H_2(\widetilde M ;\Qq/\Zz)$
  be the following composition of the duality isomorphisms and the pull back
  $$H_2(M;\Qq/\Zz)\cong  H^1(M, \partial M;\Qq/\Zz)
   \to H^1(\widetilde M, \partial \widetilde M;\Qq/\Zz) \cong  H_2(\widetilde M
   ;\Qq/\Zz).$$
 Then for any $x\in H_2(M;\Qq/\Zz)$,
 $$\ccc_{\widetilde M} (p^*(x))+1\leq n (\ccc_M(x)+1).$$
This follows from \eqref{eq1+9} and the fact that if a simple folded
surface $X$ in $M$ represents $x$, then   $p^{-1} (X)\subset
\widetilde M$ is a quasi-simple folded surface   representing
$p^*(x)$.

\subsection{Norms associated with links} A link $L$ in an oriented 3-manifold
$M$ determines a semi-norm $\Vert \cdot \Vert_{M,L}$ on $H_2(M;\Qq)$
as follows. Let $U\subset M$ be a regular neighborhood of $L$ and
$E=\overline{M-U}$ the exterior of $L$. We can embed $H_2(M;\Qq)$
into $H_2(E,\partial E;\Qq)$ via the inclusion homomorphism
$$H_2(M;\Qq)\hookrightarrow H_2(M,L;\Qq)\cong H_2(M,U;\Qq) \cong H_2(E,\partial
E;\Qq).$$ Restricting the Thurston semi-norm on $H_2(E,\partial
E;\Qq)$ to $H_2(M;\Qq)$, we obtain the semi-norm $\Vert \cdot
\Vert_{M,L}$. The arguments as above allow us to estimate the latter
semi-norm from below for compact $M$. Namely, if $L$ has $m\geq 1$
components and $h_1,...,h_m\in H=H_1(M)$ are their homology classes,
then
$$\Vert x \Vert_{M,L}\geq \spn_x (\prod_{i=1}^m (h_i-1)\, \tau)$$
for any $x\in H_2(M;\Qq)$ and  any $\tau \in Q(H)$ representing
$\tau(M)$ in the case $b_1(M)\geq 2$ and representing $[\tau] (M)$
in the case $b_1(M)=1$. A similar construction can be used to derive
a function on $H_2(M;\Qq/\Zz)$ from the function $\ccc$ on
$H_2(E,\partial E;\Qq/\Zz)$. It would be interesting to see whether
these semi-norms and   functions   may be used to  distinguish
non-isotopic links.

\subsection{Open questions} Is  the infimum in \eqref{eq1}
realizable by a simple folded surface?  Does  $\ccc $   take only
 rational values?   A positive
answer to the first question  certainly  implies a positive answer
to the second one. Similar  questions can be asked for $\CCC$.

It would be interesting
   to
compute the function $\CCC$ for the lens spaces. Is it true that for
the  lens spaces, the inequality in Theorem \ref{esb3} is   an
equality?

 \end{document}